\documentclass{article}\begin{document}
\centerline{\bf\large Solution of Polynomial Equations with Nested Radicals}\vskip .2in

\centerline{\bf Nikolaos D. Bagis}
\centerline{\bf Aristotele University of Thessaloniki (AUTH)}
\centerline{\bf Thessaloniki, Greece}

\[
\]

\centerline{\bf Abstract}

\begin{quote}
In this article we present solutions of certain  polynomial equations in periodic nested radicals. We also present a new way to solve the general tetranomial equation with new functions. As application of these new functions we solve the general sextic equation.
\end{quote}


\section{Introduction}

In (1758) Lambert considered the trinomial equation
\begin{equation}
x^m+q=x
\end{equation}
and solve it giving $x$ as a series of $q$.\\ Euler, write the same equation into, the more consistent and symmetrical form 
\begin{equation}
x^a-x^b=(a-b)\nu x^{a+b},
\end{equation}
using the transformation $x\rightarrow x^{-b}$ and setting $m=ab$, $q=(a-b)\nu$ in (1).\\
Euler gave his solution as
$$
x^n=1+n\nu+\frac{1}{2}n(n+a+b)\mu^2+\frac{1}{6}n(n+a+2b)(n+2a+b)+
$$
$$
+\frac{1}{24}n(n+a+3b)(n+2a+2b)(n+3a+b)\nu^4+
$$
\begin{equation}
+\frac{1}{120}n(n+a+4b)(n+2a+3b)(n+3a+2b)(n+4a+b)\nu^{5}+\ldots
\end{equation} 
\\
The equation of Lambert and Euler can formulated in the next (see [4] pg.306-307):\\
\\
\textbf{Theorem 1.}\\The equation 
\begin{equation}
aqx^p+x^q=1
\end{equation} 
admits root $x$ such that 
\begin{equation}
x^n=\frac{n}{q}\sum^{\infty}_{k=0}\frac{\Gamma(\{n+pk\}/q)(-q a)^k}{\Gamma(\{n+pk\}/q-k+1)k!}\textrm{, }n=1,2,3,\ldots
\end{equation}
where $\Gamma(x)$ is Euler's the Gamma function.\\
\\
 
Moreover if someone defines $\sqrt[d]{x}:=x^{1/d}$, where $d\in\textbf{Q}_{+}-\{0\}$, then the solution $x$ of (4) can given in nested radicals:   
\begin{equation}
x=\sqrt[q]{1-aq\sqrt[q/p]{1-aq\sqrt[q/p]{1-aq\sqrt[q/p]{1-\ldots}}}}.
\end{equation}
Turn now into the general quintic equation
\begin{equation}
ax^5+bx^4+cx^3+dx^2+ex+f=0. 
\end{equation}
This equation can reduced by means of a Tchirnhausen transform (see [1]) into
\begin{equation}
x^5+Ax+B=0.
\end{equation}
We can  solve this last equation using the arguments of Lambert and Euler since it is of the form (4). This can be done defining the hypergeometric function
\begin{equation}
\textbf{BR}(t)=-t\cdot{}_4F_3\left[\left\{\frac{1}{5},\frac{2}{5},\frac{3}{5},\frac{4}{5}\right\};\left\{\frac{1}{2},\frac{3}{4},\frac{5}{4}\right\};-\frac{3125 t^4}{256}\right],
\end{equation}
then according to Theorem 1:\\
\\ 
\textbf{Theorem 2.}\\ The solution of (8) is
\begin{equation}
x=\sqrt[4]{\frac{-A}{5}}\textbf{BR}\left(-\frac{\sqrt[4]{\frac{5^5}{-A^5}}}{4}B\right)
\end{equation}
and for (8) also holds
\begin{equation}
x=\sqrt[5]{-B-A\sqrt[5]{-B-A\sqrt[5]{-B-A\sqrt[5]{-B-\ldots}}}}
\end{equation}

\section{The sextic equation}

It is true that in nested radicals terminology we can do much more than that.\\ 
Consider the following general type of equation
\begin{equation}
ax^{2\mu}+bx^{\mu}+c=x^{\nu},
\end{equation}
then (12) can be written in the form
\begin{equation}
a\left(x^{\mu}+\frac{b}{2a}\right)^2-\frac{\Delta^2}{4a}=x^{\nu},
\end{equation}
where $\Delta=\sqrt{b^2-4ac}$. Hence
$$
x^{\mu}=\frac{-b}{2a}+\sqrt{\frac{\Delta^2}{4a^2}+\frac{1}{a}x^{\nu}}
$$
or
$$
x=\sqrt[\mu]{\frac{-b}{2a}+\sqrt{\frac{\Delta^2}{4a^2}+\frac{1}{a}x^{\nu}}}.
$$
Hence
\begin{equation}
x^{\nu}=\sqrt[\mu/\nu]{\frac{-b}{2a}+\sqrt{\frac{\Delta^2}{4a^2}+\frac{1}{a}x^{\nu}}}
\end{equation}
\\
From the above we get:\\
\\ 
\textbf{Theorem 3.}\\ The solution of (12) is
\begin{equation}
x=\sqrt[\mu]{\frac{-b}{2a}+\sqrt{\frac{\Delta^2}{4a^2}+\frac{1}{a}\sqrt[\mu/\nu]{\frac{-b}{2a}+\sqrt{\frac{\Delta^2}{4a^2}+\frac{1}{a}\sqrt[\mu/\nu]{\frac{-b}{2a}+\sqrt{\frac{\Delta^2}{4a^2}+...}}}}}}
\end{equation}
\\
\textbf{Corollary 1.}\\ The general sextic equation  
\begin{equation}
ax^6+bx^5+cx^4+dx^3+ex^2+fx+g=0
\end{equation}
is always solvable with periodic nested radicals.\\
\\
\textbf{Proof.}\\
We know (see [1]) that all sextic equations, of the most general form (16) are equivalent by means of a Tchirnhausen
transform 
\begin{equation}
y=kx^4+lx^3+mx^2+nx+s
\end{equation}
to the form
\begin{equation}
y^6+e_1y^2+f_1y+g_1=0
\end{equation}    
But the last equation is of the form (12) with $\mu=1$ and $\nu=6$ and have solution
\begin{equation}
y=\frac{-f_1}{2e_1}+\sqrt{\frac{\Delta_1^2}{4e_1^2}-\frac{1}{e_1}\sqrt[1/6]{\frac{-f_1}{2e_1}+\sqrt{\frac{\Delta_1^2}{4e_1^2}-\frac{1}{e_1}\sqrt[1/6]{\frac{-f_1}{2e_1}+\sqrt{\frac{\Delta_1^2}{4e_1^2}+...}}}}}
\end{equation} 
where $\Delta_1=\sqrt{f_1^2-4e_1g_1}$. Hence knowing $y$ we find $x$ from (17) and get the solvability of (16) in nested periodic radicals. More precisely it holds
$$
kx^4+lx^3+mx^2+nx+s=
$$
\begin{equation}
=\frac{-f_1}{2e_1}+\sqrt{\frac{\Delta_1^2}{4e_1^2}-\frac{1}{e_1}\sqrt[1/6]{\frac{-f_1}{2e_1}+\sqrt{\frac{\Delta_1^2}{4e_1^2}-\frac{1}{e_1}\sqrt[1/6]{\frac{-f_1}{2e_1}+\sqrt{\frac{\Delta_1^2}{4e_1^2}+...}}}}}
\end{equation}
\\
\textbf{Corollary 2.}\\ 
If $Im(\tau)>0$ and $j_{\tau}$ denotes the $j$-invariant, then
$$
\frac{1}{j_{\tau}^{5/3}}\left(R\left(e^{2\pi i \tau}\right)^{-5}-11-R\left(e^{2\pi i\tau}\right)^5\right)^{5/3}=
$$
\begin{equation}
=\sqrt[3/5]{\frac{-125}{j_{\tau}}+\sqrt{\frac{12500}{j^2_{\tau}}+\sqrt[3/5]{\frac{-125}{j_{\tau}}+\sqrt{\frac{12500}{j_{\tau}^2}+\ldots}}}},
\end{equation}
where $R(q)$ is the Rogers-Ramanujan continued fraction:
\begin{equation}
R(q)=\frac{q^{1/5}}{1+}\frac{q}{1+}\frac{q^2}{1+}\frac{q^3}{1+}\ldots\textrm{, }|q|<1
\end{equation}
\\
\textbf{Proof.}\\
Equation (12) for $a=1/j_{\tau}^{1/3}$, $b=250/j_{\tau}^{1/3}$, $c=3125/j_{\tau}^{1/3}$ and $\mu=1$, $\nu=5/3$ takes the form 
\begin{equation}
x^{2}+250x+3125=j_r^{1/3}x^{5/3},
\end{equation}
which is a simplified form of Klein's equation for the icosahedron (see [2] and [3]) and have solution $x=Y_{\tau}=R(q^2)^{-5}-11-R(q^2)^5$. From Theorem 3 we can express $Y_{\tau}$ in nested-periodical radicals. This completes the proof.\\  
\\
\textbf{Note.} For more details about equation $ax^2+bx+\frac{b^2}{20a}=C_1x^{5/3}$ one can see [2].

\section{A solution of a tetranomial polynomial}

Consider now the generalization of the Bring radical function as follows:\\
Set $p=\nu$ and $q=\nu-1$, $x=X^{-1}$ in (4), then it becomes
$$
a(\nu-1)X^{-\nu}+X^{-\nu+1}=1
$$
Multiply both sides with $X^{\nu}$ and set $a\rightarrow a/(\nu-1)$, then 
\begin{equation}
X^{\nu}=X+a
\end{equation}
This last equation according to Theorem 1 have solution
\begin{equation}
X=\left(\frac{1}{\nu-1}\sum^{\infty}_{n=0}\frac{\Gamma[(1+\nu n)/(\nu-1)]}{\Gamma[(-1+\nu n)/(\nu-1)-n+1]}\frac{(-a)^n}{n!}\right)^{-1},
\end{equation}
which is a hypergeometric function that we call $\nu-$th order Bring radical. Also the power $X^{\mu}$ has expansion
\begin{equation}
X^{\mu}=\left(\frac{\mu}{\nu-1}\sum^{\infty}_{n=0}\frac{\Gamma[(\mu+\nu n)/(\nu-1)]}{\Gamma[(\mu+\nu n)/(\nu-1)-n+1]}\frac{(-a)^n}{n!}\right)^{-1}.
\end{equation}
A simple way to write (25),(26) is
\begin{equation}
X=\textbf{BR}_{\nu}(a)=\left(\sum^{\infty}_{n=0}\left(\begin{array}{cc}
	\frac{\nu n+1}{\nu-1}\\
	n
\end{array}
\right)\frac{(-a)^n}{\nu n+1}\right)^{-1}
\end{equation}
and from Theorem 1: 
\begin{equation}
X^{\mu}=\lambda_1\sum^{\infty}_{n=0}\frac{(-1)^n}{k_1n+\lambda_1}\left(
\begin{array}{cc}
	k_1n+\lambda_1\\
	n
\end{array}
\right)a^n,
\end{equation}
where $k_1=\frac{\nu}{\nu-1}$ and $\lambda_1=\frac{-\mu}{\nu-1}$. Hence we lead to the following\\
\\
\textbf{Definition 1.}\\
We call ''General Bring Radical'' the function
\begin{equation}
\textbf{BR}_{\nu,\mu}(a):=\left(\textbf{BR}_{\nu}(a)\right)^{\mu}=\lambda_1\sum^{\infty}_{n=0}\frac{(-1)^n}{k_1n+\lambda_1}\left(\begin{array}{cc}
	k_1n+\lambda_1\\
	n
\end{array}\right)a^n,
\end{equation}
where $k_1=\frac{\nu}{\nu-1}$ and $\lambda_1=\frac{-\mu}{\nu-1}$.\\
\\

Assume now the following equation
\begin{equation}
kX^{\mu}=X^{\nu}-X-l.
\end{equation}
If $H_{\nu}(u)$ is a function such that
\begin{equation}
H_{\nu}(u)^{\nu}-H_{\nu}(u)-l=u,
\end{equation} 
then equation (30) have solution $X_2=H_{\nu}(u_0)$ such that 
$$
H_{\nu}\left(u_0\right)^{\nu}-H_{\nu}\left(u_0\right)-l=u_0=kH_{\nu}\left(u_0\right)^{\mu}.
$$
But 
\begin{equation}
H_{\nu}(u)=\textbf{BR}_{\nu}(l+u)\textrm{, }\forall u\in D,
\end{equation}
where $D$ is a suitable domain. Hence
$$
X_2=H_{\nu}(u_0)=\textbf{BR}_{\nu}\left(l+u_0\right)=
$$
\begin{equation}
=\textbf{BR}_{\nu}\left(l+k\textbf{BR}_{\nu,\mu}\left(l+k\textbf{BR}_{\nu,\mu}(l+\ldots)\right)\right)
\end{equation}
and we get the following\\
\\
\textbf{Theorem 4.}\\
The equation
\begin{equation}
X^{\nu}=kX^{\mu}+X+l
\end{equation}
admits solution $X_2$ given by
\begin{equation}
\xi=l+k\textbf{BR}_{\nu}(\xi)^{\mu}\textrm{, }X_2=\textbf{BR}_{\nu}(\xi)
\end{equation}
Setting this into nested Bring radicals we have
\begin{equation}
X_2=\textbf{BR}_{\nu}\left(l+k\textbf{BR}_{\nu,\mu}\left(l+k\textbf{BR}_{\nu,\mu}(l+\ldots)\right)\right)
\end{equation}
\\
\textbf{Application.}\\
Assume the equation
\begin{equation}
X^{\mu}=X^2-X-1.
\end{equation}
Then $\mu\in\textbf{R}$ and $\nu=2$ and if 
\begin{equation}
f_0(x):=\textbf{BR}_{2,1}(x)=\frac{1+\sqrt{1+4x}}{2}
\end{equation}
and
\begin{equation}
f(x):=\textbf{BR}_{2,\mu}(x)=\left(\frac{1+\sqrt{1+4x}}{2}\right)^{\mu}
\end{equation}
the above equation (37) admits solution 
\begin{equation}
X_2=f_0\left(1+f\left(1+f(1+\ldots)\right)\right).
\end{equation}
\\
\textbf{Theorem 5.}\\
Assume the equation
\begin{equation}
X^{\nu}=kX^{\mu}+X+l,
\end{equation}
with $\nu>\mu$ and $|k|,|l|,|kl^{-1}|<1$. Set $k_1=\frac{\nu}{\nu-1}$, $\lambda_1=\frac{-\mu}{\nu-1}$,
\begin{equation}
\textbf{f}(s;x):=\sum^{\infty}_{n=1}\left(\begin{array}{cc}
	s\\
	n
\end{array}\right)\left(\begin{array}{cc}
	k_1s+\lambda_1n\\
	s
\end{array}\right)\frac{n}{(s-n+1)(k_1 s+\lambda_1n)}x^n
\end{equation}
and
\begin{equation}
\textbf{g}(w):=w+w\lambda_1\sum^{\infty}_{s=0}(-1)^s\textbf{f}\left(s;kw^{-1}\right)w^s.
\end{equation}
Then a solution of (41) is
\begin{equation}
X_2=\textbf{BR}_{\nu}\left(\textbf{g}(l)\right).
\end{equation}
\\
\textbf{Proof.}\\
It is clear that if $X_2$ is the desired root of the equation (41), then also $X_2=\textbf{BR}_{\nu}\left(\xi\right)$, where $\xi=l+k\textbf{BR}_{\nu,\mu}\left(\xi\right)$. But according to Lagrange theorem (see [6] pg 133), we have
\begin{equation}
\xi=l+\sum^{\infty}_{n=1}\frac{k^n}{n!}\frac{d^{n-1}}{dl^{n-1}}\left(\textbf{BR}_{\nu,\mu}(l)\right)^n
\end{equation}
and  $\left(\textbf{BR}_{\nu,\mu}(x)\right)^{n}=\textbf{BR}_{\nu,\mu n}(x)$. Using expansion (29) we arrive to a double sum we rearrange and get the result.\\
\\
\textbf{Example.}\\
Assume that $\nu=2$, $\mu=1/2$ and $l=1/5$, $k=1/5^2$, then 
\begin{equation}
X^{1/2}=-5-25X+25X^2
\end{equation} 
and 
\begin{equation}
\textbf{BR}_{2,1}(x)=\frac{1+\sqrt{1+4x}}{2}
\end{equation}
The polynomial $\textbf{f}\left(s;kl^{-1}\right)$ is
$$
\textbf{f}\left(s;\frac{1}{5}\right)=\frac{s C_{2 s-1,s}}{50 s-25}\cdot{ }_3F_2\left(1-s,\frac{1}{2}-\frac{s}{2},1-\frac{s}{2};\frac{3}{2},2-2 s;\frac{1}{25}\right)+
$$
\begin{equation}
+\frac{2 C_{2s-\frac{1}{2},s}}{20 s-5}\cdot{ } _3F_2\left(\frac{1}{2}-s,\frac{1}{2}-\frac{s}{2},-\frac{s}{2};\frac{1}{2},\frac{3}{2}-2 s;\frac{1}{25}\right)
\end{equation}
and $\textbf{g}$ will be
\begin{equation}
\textbf{g}(l)=\frac{1}{5}-\frac{1}{10}\sum^{\infty}_{s=0}(-1)^s5^{-s}\textbf{f}\left(s;\frac{1}{5}\right),
\end{equation}
where we have denoted
\begin{equation}
C_{n,m}:=\left(\begin{array}{cc}
	n\\
	m
\end{array}\right).
\end{equation}
Hence a solution of (46) will be
\begin{equation}
X_2=\textbf{BR}_{2,1}\left(\textbf{g}(l)\right).
\end{equation}
\\
\textbf{Corallary 3.}\\
1) Every sextic equation can be set in the form
\begin{equation}
x^6+ax^2+bx+c=0
\end{equation}
2) Under the change of variable $x=\sqrt[5]{b}X$ the equation (52) takes the form
\begin{equation}
X^6+ab^{-4/5}X^2+X+cb^{-6/5}=0
\end{equation}
3) If $k=-ab^{-4/5}$, $l=-cb^{-6/5}$, $|kl^{-1}|=|ab^{2/5}c^{-1}|<1$ and $|k|\leq\frac{1}{4}$, $|l|\leq\frac{1}{3}$, then equation (53) is solvable with the functions (42),(43),(44) introduced in Theorem 5.\\
\\
\textbf{Theorem 6.}\\
Assume the equation 
\begin{equation}
X^{\nu}=mX^{\lambda}+kX^{\mu}+X+l,
\end{equation}
where $\nu>\mu$ and $|k|,|l|,|kl^{-1}|<1$. Assume also the equation 
\begin{equation}
X^{\nu}=kX^{\mu}+X+l.
\end{equation}
Then a solution of (54) is 
\begin{equation}
X_3=\textbf{BR}_{\nu}(\textbf{g}(l+m\textbf{BR}_{\nu,\lambda}(\textbf{g}(l+m\textbf{BR}_{\nu,\lambda}(\textbf{g}(l+\ldots)))))),
\end{equation}
with ''denested'' equation 
\begin{equation}
\xi=l+m\textbf{BR}_{\nu,\lambda}\left(\textbf{g}(\xi)\right)\textrm{ and }X_3=\textbf{BR}_{\nu}\left(\textbf{g}(\xi)\right).
\end{equation}
Here the functions $\textbf{g}$ and $\textbf{f}$ are exactly that of Theorem 5.\\ 
\\
\textbf{Proof.}\\
Set $H(u)$ be a function such that 
\begin{equation}
-kH(u)^{\mu}+H(u)^{\nu}-H(u)-l=u.
\end{equation}
The function $H(u)$ can determined using Theorem 5.  
Assume also that equation (54) admits solution $X_3=H(u_0)$ and for this $u=u_0$ we have 
\begin{equation}
mH(u_0)^{\lambda}=-kH(u_0)^{\mu}+H(u_0)^{\nu}-H(u_0)-l=u_0.
\end{equation}
Then
$$
H(u)=\textbf{BR}_{\nu}(\textbf{g}(l+u))
$$
and hence 
$$
H(u_0)=\textbf{BR}_{\nu}(\textbf{g}(l+u_0))\textrm{, with }
mH(u_0)^{\lambda}=u_0.
$$
Hence a root of (54) is
\begin{equation}
X_3=H(u_0)=\textbf{BR}_{\nu}(\textbf{g}(l+m\textbf{BR}_{\nu,\lambda}(\textbf{g}(l+m\textbf{BR}_{\nu,\lambda}(\textbf{g}(l+\ldots)))))).
\end{equation}
\\
\textbf{Theorem 7.}\\
According to Lagrange inversion theorem equation (54) admits root
\begin{equation}
X_3=\textbf{BR}_{\nu}\left(\textbf{g}\left(l+m\sum^{\infty}_{n=0}\frac{m^n}{(n+1)!}\left[\frac{d^n}{dh^n}\textbf{BR}_{\nu,n\lambda}(\textbf{g}(h))\right]_{h=l}\right)\right)
\end{equation}

\newpage

\centerline{\bf References}\vskip .2in

\noindent

[1]: J.V. Armitage and W.F. Eberlein. 'Elliptic Functions'. Cambridge University Press. (2006)

[2]: N.D. Bagis. 'On a general sextic equation solved by the Rogers-Ramanujan continued fraction'. arXiv:1111.6023v2 [math.GM]. (2012)

[3]: W. Duke. 'Continued fractions and Modular functions'. Bull. Amer. Math. Soc. (N.S.), 42. (2005), 137-162.

[4]: B.C. Berndt. 'Ramanujan`s Notebooks Part I'. Springer Verlag, New York. (1985).

[5]: M.L. Glasser. 'The Quadratic Formula Made Hard or A Less Radical Approach to Solving Equations'. arXiv:math/94112244v1 [math.CA] 16 Nov (1994).

[6]: E.T. Whittaker and G.N. Watson. 'A course on Modern Analysis'. Cambridge U.P, (1927).

\end{document}